\documentclass[10pt]{article}
\usepackage[dvipdfmx]{graphicx}
\def\qed{\hfill \hbox{${\vcenter{\vbox{
   \hrule height 0.4pt\hbox{\vrule width 0.4pt height 6pt
   \kern5pt\vrule width 0.4pt}\hrule height 0.4pt}}}$}}
\newcommand{\maru}[1]{{\ooalign{\hfil#1\/\hfil\crcr
   \raise.167ex\hbox{\mathhexbox20D}}}}
\newcommand{\bysame}{%
   \leavevmode\hbox to 6em{\hrulefill}\,}

\begin{document}

\centerline{\large\bf On composite twisted torus knots}
\vskip 2mm 

\centerline{by} 
\vskip 2mm 

\centerline{\bf Kanji Morimoto} 
\vskip 5mm 

\centerline{\large\it Dedicated to Professor Akio Kawauchi for his 60th birthday} 
\vskip 5mm 

\centerline{Department of IS and Mathematics, Konan University} 
\vskip -1mm 

\centerline{Okamoto 8-9-1, Higashi-Nada, Kobe 658-8501, Japan}
\vskip -1mm 

\centerline{morimoto@konan-u.ac.jp} 
\vskip 5mm 

{\bf Abstract.} \ In the present note, we will show that there are infinitely many composite twisted torus knots. 
\vskip 5mm 

Keywords and phrases \ : \ twisted torus knots, composite, prime 
 
2000 Mathematics Subject Classification \ : \ 57M25, 57N10 
\vskip 10mm 

{\bf 1. Introduction} 
\vskip 3mm 

\quad Let $K$ be a knot in the 3-sphere $S^3$. Suppose $K$ is the connected sum of two non-trivial knots $K_1$ and $K_2$. Then we say that $K$ is a composite knot, and denote it by $K = K_1 \# K_2$. Otherwise we say that $K$ is a prime knot. 

\quad Let $p, q, r, s$ be integers such that $p > r > 1$, $q > 0$, gcd$(p, q)=1$ and let $T(p, q)$ be the torus knot of type $(p, q)$ in $S^3$. For the definition of torus knots $T(p, q)$ we refer to [7]. Add $s$ times full twists on mutually parallel $r$ strands in $T(p, q)$. Then according as [1], we call the knot obtained by this operation a twisted torus knot of type $(p, q ; r, s)$ and denote it by $T(p, q ; r, s)$ as illustrated in Figure 1.  

\begin{figure}[htbp]
\hskip 23mm 
\includegraphics[width=6.5cm]{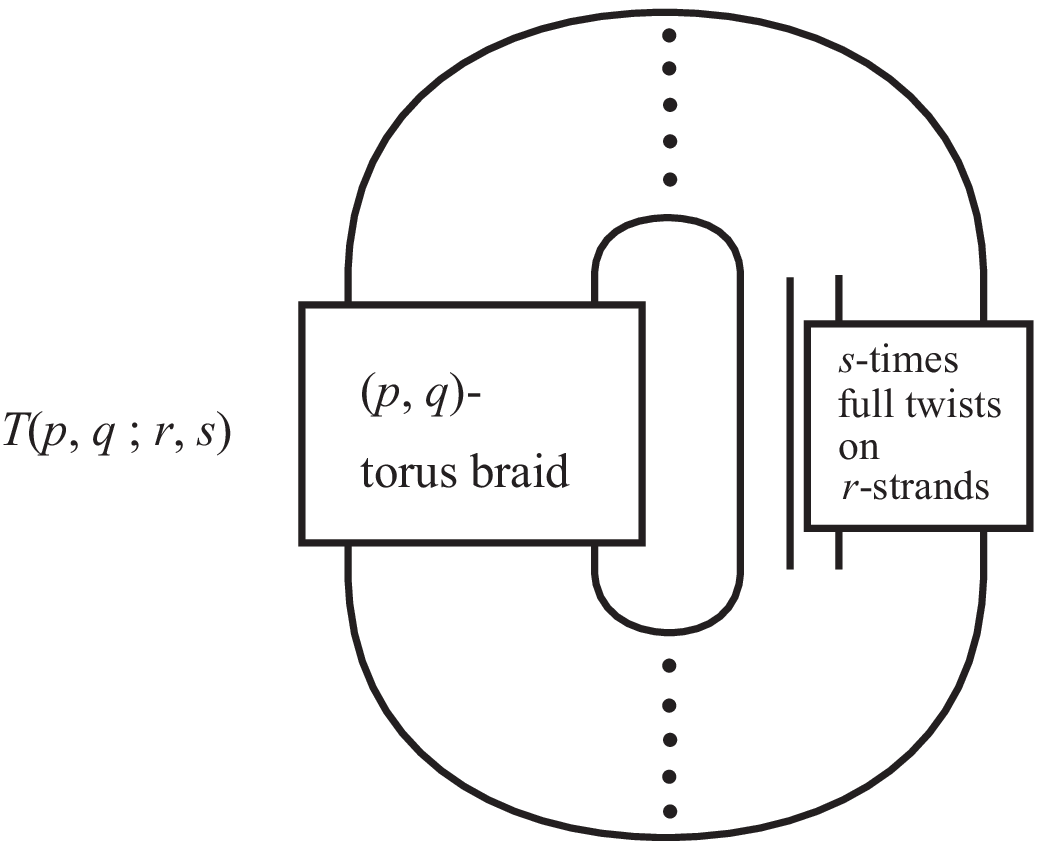}

\center{Figure 1}
\end{figure}

\quad Twisted torus knots are deeply related to unexpected Dehn surgeries. In fact, the famous hyperbolic pretzel knot $P(-2,3,7)$ found by Fintushel--Stern in [2] is the twisted torus knot $T(5, 3 ; 2, 1)$. In addition, many similar hyperbolic twisted torus knots have been found in [1]. Moreover, twisted torus knots have interesting properties in the additivity of tunnel numbers of knots as in [4]. Therefore, the family of twisted torus knots has been considered as an important class in studying of knot theory.  

\quad By a little observation, we see that $T(p, q ; 2, s)$ has tunnel number one for any choice of $p, q, s$, and is prime by [6]. Moreover, J. H. Lee has recently shown in [3] that $T(p, q ; 3, s)$ has also tunnel number one for any choice of $p, q, s$, and is prime by [6] again. In fact, so far, no composite twisted torus knot has been known. Therefore, we need to ask if there are composite twisted torus knots. In the present note, we will answer to this question as follows : 
\vskip 3mm 

{\bf Theorem 1.} \ \it Suppose $p=(a+1)(k_1+k_2)+1$, $q=a(k_1+k_2)+1$, $r=p-k_1$ and $s=-1$ for some integers $a>0$, $k_1>1$ and $k_2>1$. Then $T(p, q ; r, s)$ is the connected sum of $T(k_1, ak_1+1)$ and $T(k_2, -(a+1)k_2-1)$. \rm 
\vskip 3mm 

{\bf Examples.} \ $(1)$ Put $a=1$, $k_1=k_2=2$, then $p=9$, $q=5$, $r=7$ and $T(9, 5 ; 7, -1) \cong T(2, 3) \# T(2, -5)$. 

$(2)$ Put $a=2$, $k_1=4, k_2=2$, then $p=19$, $q=13$, $r=15$ and $T(19, 13 ; 15, -1) \cong T(4, 9) \# T(2, -7)$.   
\vskip 3mm 

\quad By the way, it is well known that composite knots have essential tori in the exteriors. Concerning the conditin for twisted torus knots to have essential tori in the exteriors, we have shown in [5] that for any composite number $r = km \ (k>1, m>1)$, there are infinitely many twisted torus knots $T(p, q ; r, s)$ which have essential tori in the exteriors. Moreover, we have shown that those knots are cable knots along some torus knots, and are prime. Therefore, by these results, we need to consider the following problem : 
\vskip 3mm 

{\bf Problems.} \ \it $(1)$ Characterize the knot types of composite twisted torus knots. In particular, we conjecture that the condition in Theorem $1$ is also a necessary condition for twisted torus knots to be composite knots. 

$(2)$ Characterize the knot types of prime twisted torus knots with essential tori. \rm 
\vskip 3mm 
  
{\bf 2. Proof of Theorem 1} 
\vskip 3mm 
 
\quad Let $K=T(p, q ; r, s)$ be the knot as in Theorem 1, i.e., $p=(a+1)(k_1+k_2)+1$, $q=a(k_1+k_2)+1$, $r=p-k_1$ and $s=-1$ for some integers $a>0$, $k_1>1$ and $k_2>1$. Then we can regard $p, q, r$ as follows : 

\quad $p=(a+1)(k_1+k_2)+1=k_1+k_2+k_1+k_2+ \cdots +k_1+(k_2+1)$, 

\quad $q=a(k_1+k_2)+1=k_1+k_2+ \cdots +k_1+(k_2+1)$, 

\quad $r=p-k_1=k_1+k_2+k_1+ \cdots +k_2+(k_2+1)$. 

\quad Then we can divide $p$ strings into $a+1$ bunchs of $k_1$ strings, $a$ bunchs of $k_2$ strings and one bunch of $k_2+1$ strings, can divide $q$ strings into $a$ bunchs of $k_1$ strings, $a-1$ bunchs of $k_2$ strings and one bunch of $k_2+1$ strings and can divide $r$ strings into $a$ bunchs of $k_2$ strings, $a$ bunchs of $k_1$ strings and one bunch of $k_2+1$ strings as in Figure 2, where Figure 2 is the case of $a=2, k_1=4, k_2=2$, i.e., $K=T(19, 13 ; 15, -1)$. 

\quad First, deform the first bunch of $k_1$ strings in the $p$ strings as in Figure 3(1), and then deform the second bunch, the third bunch, $\cdots$, the $a$th bunch of $k_1$ strings, and finally decompose the knot $T(p, q ; r, s)$ into two knots at the place indicated in Figure 3(2). 

\quad Next, take the knot which consists of $k_1$ strings from the two knots obtained in Figure 3(2) as in Figure 4(1). Then we can see that the knot in Figure 4(1) is a torus knot which consists of $k_1$ strings with $a$ times full twists and $\displaystyle {2\pi \over {k_1}}$ rotation as in Figure 4(2). Thus we see that the knot is the torus knot of type $(k_1, ak_1+1)$. 

\quad Finally, take the other knot in Figure 3(2). Then we can see that the knot in Figure 5(1) is a torus knot which consists of $(a+1)k_2+1$ strings with $\displaystyle {2\pi \over {(a+1)k_2+1}}\cdot(-k_2)$ rotation as in Figure 5(2) because $ak_2+1-((a+1)k_2+1) = -k_2$. Thus we see that the knot is the torus knot of type $((a+1)k_2+1, -k_2)$. In addition, this knot is the same knot as the torus knot of type $(k_2, -(a+1)k_2-1)$. This completes the proof of Theorem1. In the case of Figure 2, we have $K=T(19, 13 ; 15, -1) = T(4, 9) \# T(7, -2) = T(4, 9) \# T(2, -7)$. \qed

\vskip 5mm 

\begin{figure}[htbp]
\hskip 23mm 
\includegraphics[width=8cm]{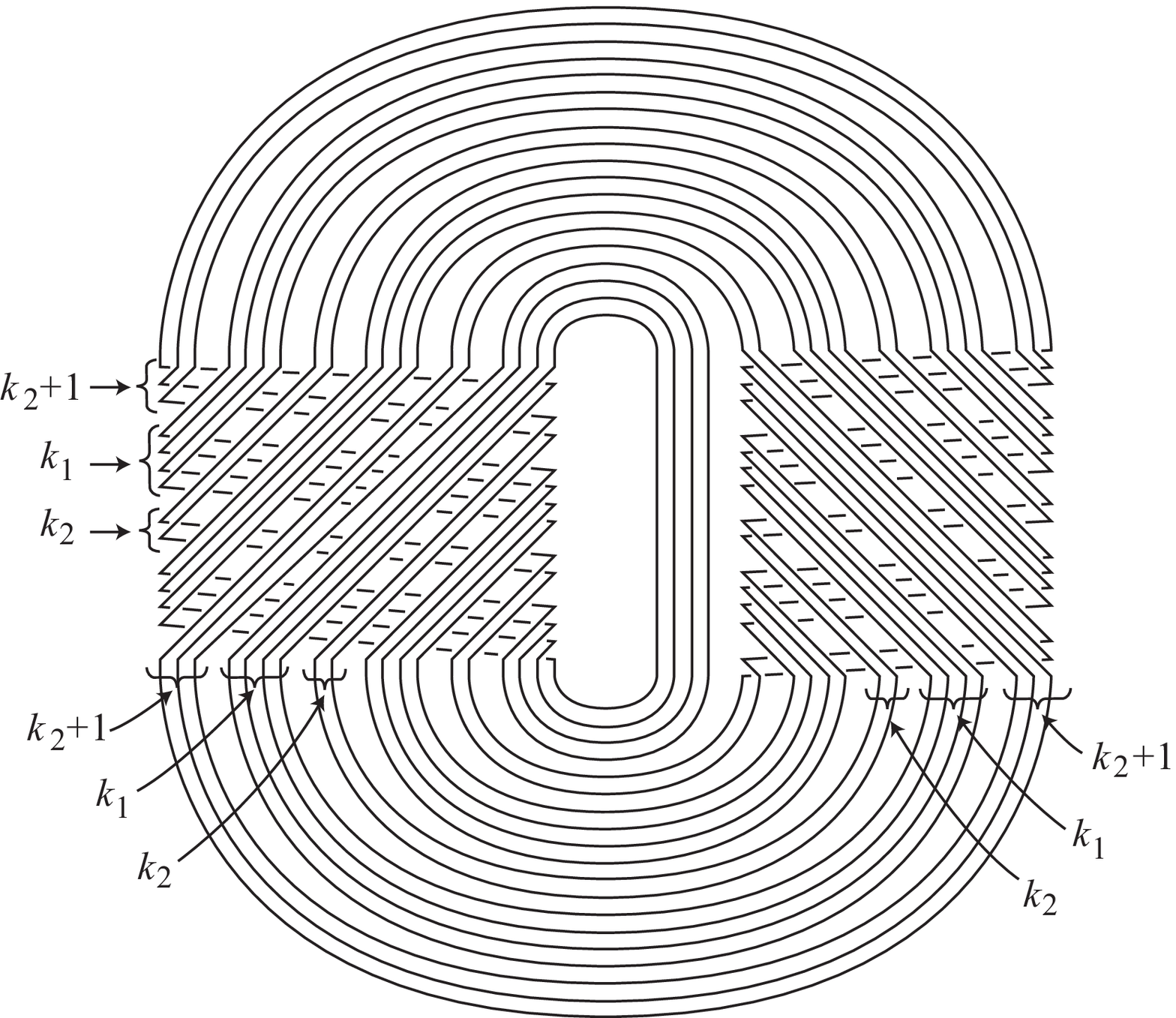}

\center{Figure 2}
\end{figure}

\begin{figure}[htbp]
\hskip 10mm 
\includegraphics[width=12.5cm]{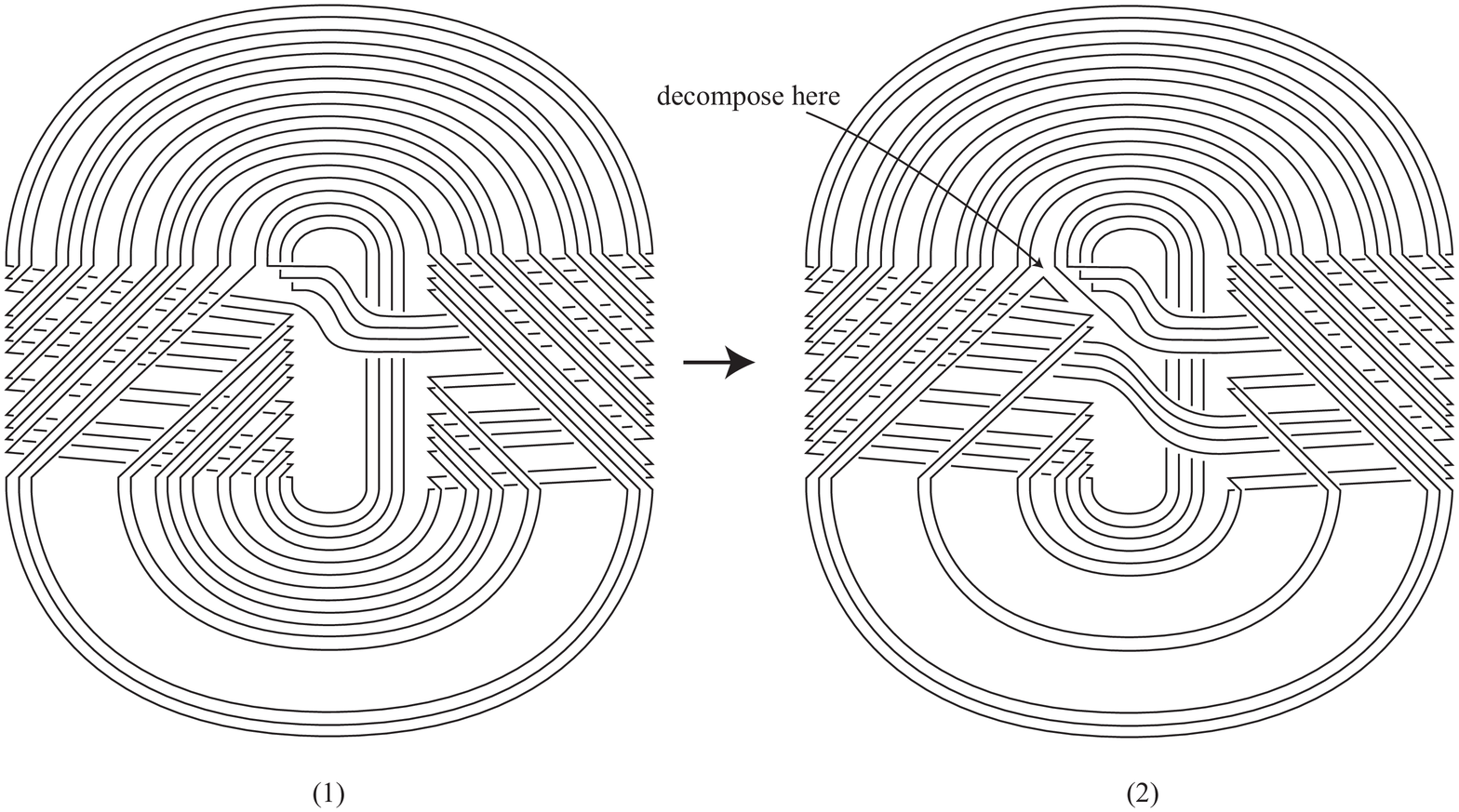}

\center{Figure 3}
\end{figure}

\begin{figure}[htbp]
\hskip 15mm 
\includegraphics[width=11cm]{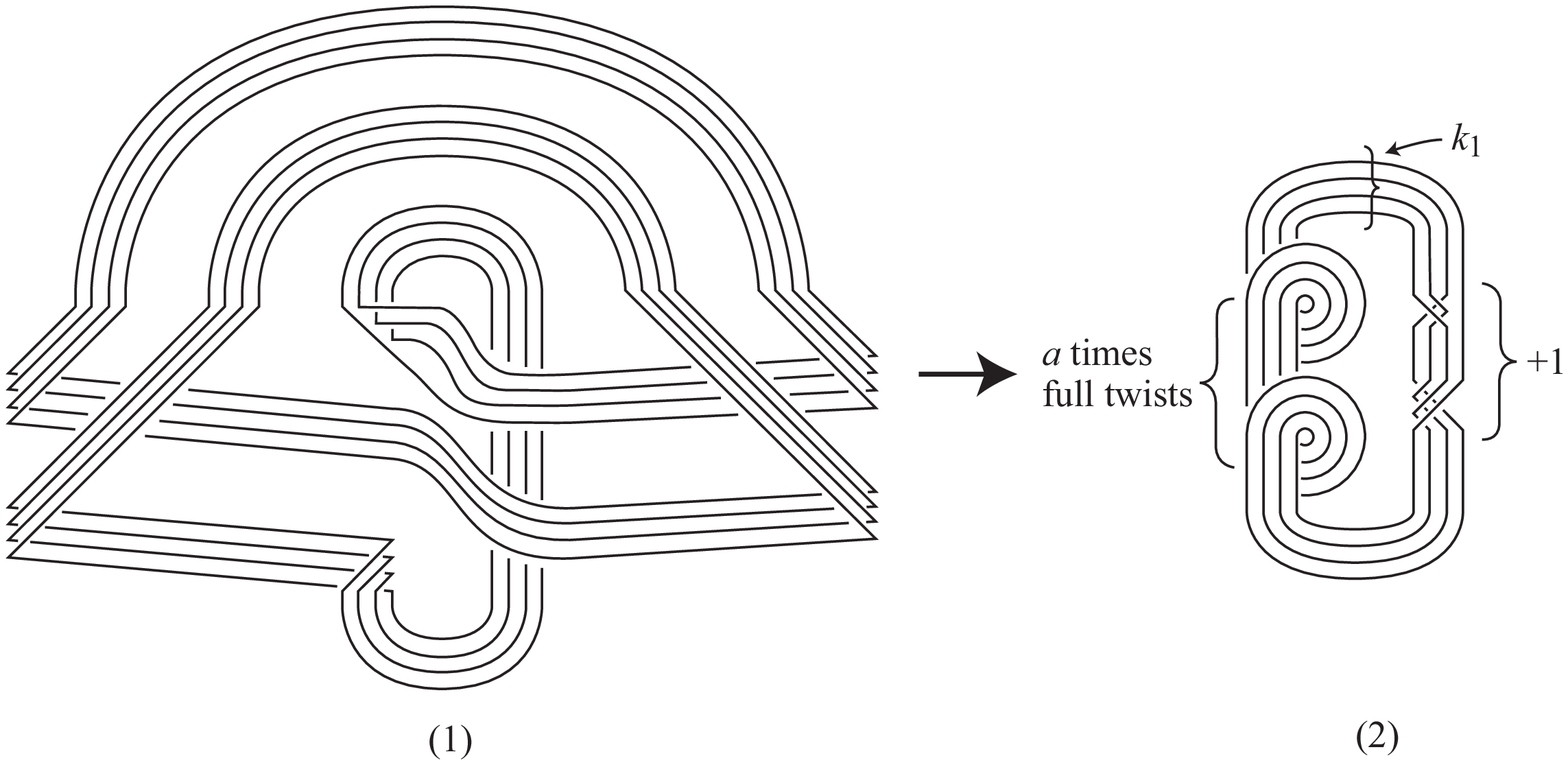}

\center {Figure 4}
\end{figure}

\begin{figure}[htbp]
\hskip 15mm 
\includegraphics[width=11cm]{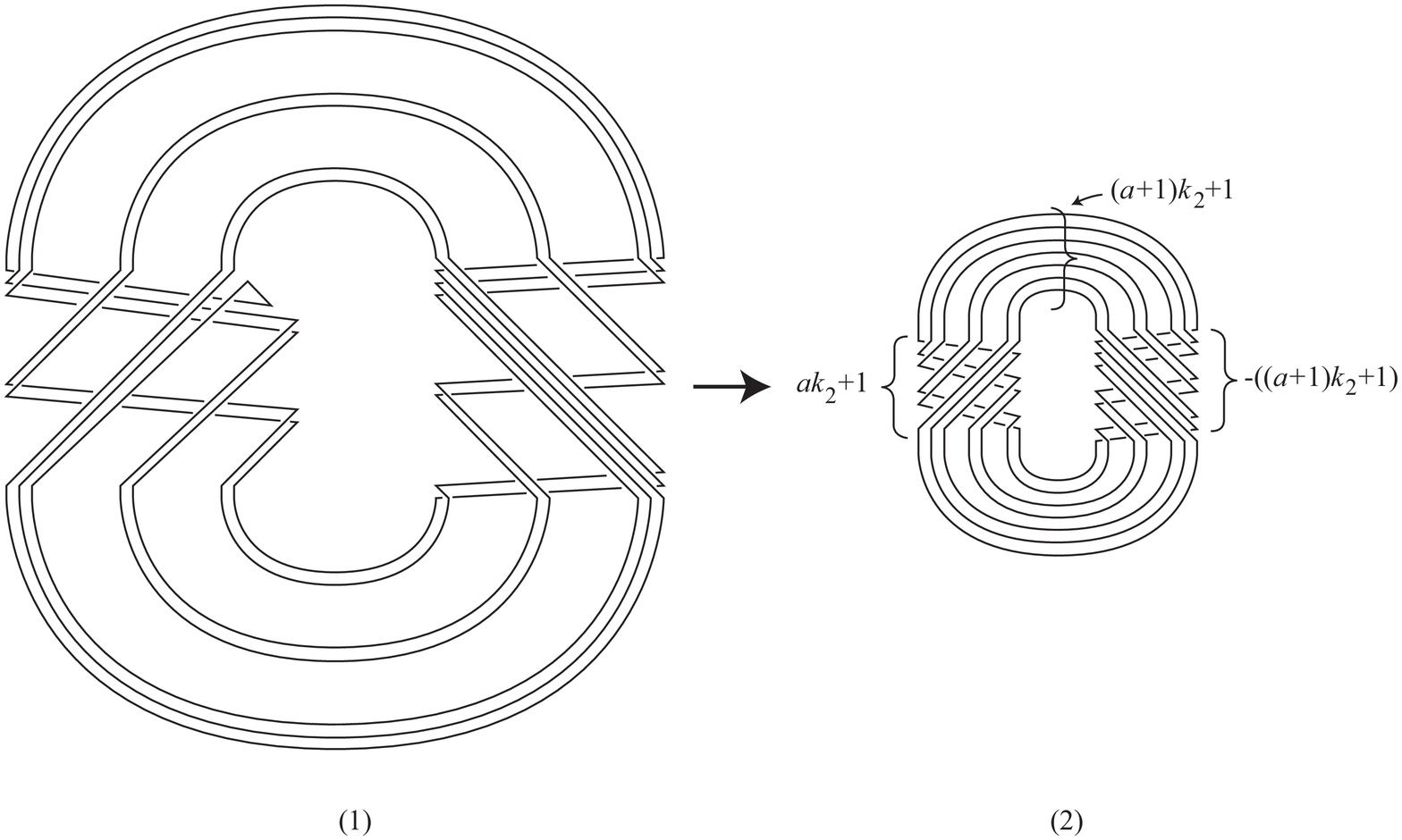}

\center {Figure 5}
\end{figure}

\vfill\eject

{\bf References} 

[1] \ P. J. Callahan, J. C. Dean and J. R.Weeks,
{\it The simplest hyperbolic knots}, 
\vskip -1mm \hskip 6mm 
J. Knot Theory Ramifications, {\bf 8} (1999) 279--297. 
\vskip -1mm 
[2] \ R. Fintushel and R. Stern,
{\it Constructing lens spaces by surgery on knots},
\vskip -1mm \hskip 6mm 
Math. Z. {\bf 175} (1980) 33--51. 
\vskip -1mm 
[3] \ J. H. Lee, \ {\it Twisted torus knots $T(p, q ; 3, s)$ are tunnel number one}, preprint. 
\vskip -1mm 
[4] \ K. Morimoto, M. Sakuma and Y. Yokota, \ {\it Examples of tunnel number one} 
\vskip -1mm \hskip 6mm
{\it knots which have the property `` $1 + 1 = 3$ ''}, \ Math. Proc. Camb. Phil. Soc. 
\vskip -1mm \hskip 6mm 
{\bf 119} (1996)  113--118 
\vskip -1mm 
[5] \ K. Morimoto and Y. Yamada, \ {\it A note on essential tori in the exteriors of torus} 
\vskip -1mm \hskip 6mm 
{\it knots with twists}, Kobe J. Math., {\bf 26} (2009) 29--34. 
\vskip -1mm 
[6] \ F. H. Norwood, \ {\it Every two generator knot is prime}, Proc. A. M. S., {\bf 86} (1982) 
\vskip -1mm \hskip 6mm 
143-147. 
\vskip -1mm 
[7] \ D. Rolfsen, \ {\it Knots and Links}, AMS Chelsea Publishing (2003)

\end{document}